\documentclass[11pt]{amsart}
\usepackage{amsmath}
\usepackage{amssymb}
\usepackage{amscd}
\newtheorem{theo}{Theorem}
\newtheorem{df}{Definition}
\newtheorem{cor}{Corollary}

\renewcommand{\Re}{\operatorname{Re}}
\renewcommand{\Im}{\operatorname{Im}}
\newcommand{\id}{\operatorname{id}}
\newcommand{\C}{{\mathbb C}}
\newcommand{\R}{{\mathbb R}}

\newcommand{\I}{\operatorname{i}}

\newcommand{\Ann}{\operatorname{Ann}}
\newcommand{\Ad}{\operatorname{Ad}}

\def\sideremark#1{\ifvmode\leavevmode\fi\vadjust{%            The remark
\vbox to0pt{\hbox to 0pt{\hskip\hsize\hskip1em%               will appear only
\vbox{\hsize3cm\tiny\raggedright\pretolerance10000%          on the side
\noindent #1\hfill}\hss}\vbox to8pt{\vfil}\vss}}}%           in 3cm

\begin{document}
\title[Cartan connection on Engel CR manifolds]{Canonical
Cartan connection and holomorphic invariants on Engel CR
manifolds}
\author{Valeri\u{i} Beloshapka \and  Vladimir Ezhov \and Gerd Schmalz}

\address{(V. Beloshapka) Moscow State University, Dept. of Mechanics
and Mathematics, Vorobyovy Gory, 119899 Moscow, Russia}
\email{vkb@strogino.ru}
\address{(V. Ezhov) School of Mathematics and Statistics,
University of South Australia, Mawson Lakes Boulevard, Mawson
Lakes, SA 5095, Australia} \email{vladimir.ejov@unisa.edu.au}
\address{(G. Schmalz) School of Mathematics, Statistics and Computer
Science, University of New England, Armidale, NSW 2351,
Australia// Uniwersytet Warmi\'nsko-Mazurski, Wydzia{\l}
Matematyki-Informatyki, ul. \.Zo{\l}nierska 14, 10-561 Olsztyn,
Poland}\email{gerd@turing.une.edu.au}
\thanks{This work was carried out in the framework of Australian Research
Council Discovery Project DP0450725}
\keywords{CR manifolds, Engel manifolds, Cartan connection}
\subjclass[2000]{Primary 32V40, 53B12 ; Secondary 34C14}

\begin{abstract}
We describe a complete system of invariants for 4-dimensional CR
manifolds of CR dimension 1 and codimension 2 with Engel CR
distribution by constructing an explicit canonical Cartan
connection. We also investigate the relation between the Cartan
connection and the normal form of the defining equation of an
embedded Engel CR manifold.

\end{abstract}
\maketitle
\setcounter{section}{-1}
\section{Introduction}
The equivalence problem for real hypersurfaces in 2-dimensional
complex space with respect to holomorphic mappings became one of
the first impressive applications of E. Cartan method of moving
frames. Cartan solved this problem in 1932 (see \cite{Car32}). In
modern terminology his approach can be described as follows: Let
$M$ be a manifold with some additional geometric data (in our case
the data come from the embedding of a real manifold into a complex
manifold). Construct a canonical principal bundle $\mathcal G$
over $M$ and a distinguished frame on $\mathcal G$ which depend
only on the specified geometric data. This means that for two
isomorphic structures $M_1$ and $M_2$ a mapping that establishes
the isomorphy lifts to a unique mapping of the canonical bundles
$\mathcal G_1$ and $\mathcal G_2$ that maps the distinguished
frames at $\mathcal G_1$ to the distinguished frames at $\mathcal
G_2$. Hence, the initial equivalence problem reduces to an
equivalence problem for manifolds with distinguished frames. Such
structure is called {\em absolute parallelism} and its geometry is
well-understood (see, e.g., \cite{Kob}). Cartan connections were constructed for nondegenerate CR-manifolds of hypersurface type (see \cite{Tan}, \cite{CM}), for 2-codimensional CR-manifolds of CR-dimension $2$ of hyperbolic and elliptic types (\cite{EIS}, \cite{SSl}, \cite{SSl1} \cite{CSa}, \cite{SSp}). For 5-dimensional uniformly degenerate CR-manifolds of hypersurface type a parallelism (with respect to a subgroup of the structure group) was obtained in \cite{Eb} (see also \cite{Eber}). 

A Cartan connection is a pair of a canonical bundle and a
distinguished frame field that satisfies certain invariance
conditions (see Section \ref{prelim} for the precise definition).
Those conditions are very useful in the study of automorphisms of
manifolds with given geometric structure.

In this paper we will describe a Cartan connection for some class
of CR manifolds that are related to third order ODE.

\begin{df}
A CR manifold is a smooth manifold $M$ equipped with a
distribution $D\subset TM$ of even rank $2n$ and a smooth field
$J$ of endomorphisms $J_a: D_a \to D_a$ with $J_a^2=-\id$. A CR
manifold is supposed to satisfy the following integrability
condition: The complex distribution spanned by $X-\I JX$ for local
sections $X$ of $D$ is involutive.

The number $n$ is called CR dimension of $M$ and $k=\dim M-2n$ is
called the codimension of $M$.
\end{df}

CR manifolds appear naturally as embedded real submanifolds of
complex manifolds. The distribution $D$ is then defined as
$D=TM\cap \I TM$ and $J$ is the restriction of the complex
structure in the ambient complex manifold to $D$.

S.S. Chern and J. Moser \cite{CM} and N. Tanaka \cite{Tan}
constructed canonical Cartan connections for Levi-nondegenerate
1-codimensional CR manifolds of arbitrary dimension. Similar
results for so-called hyperbolic and elliptic CR manifolds of
codimension two in $\C^4$ were obtained by G. Schmalz and J.
Slov\'ak \cite{SSl, SSl1} (see also \cite{EIS,CSa,SSp}). These
geometries fit into the general concept of {\em parabolic
geometry} introduced by Ch. Fefferman \cite{Fef} (see also
\cite{Yama}, \cite{CS}). This is due to the following features
they share:

\begin{enumerate}
\item {\bf Nondegeneracy.} The algebraic bracket (see Subsection \ref{brac})
$$\mathcal L_a: D_a \otimes D_a \to T_a/D_a$$
that assigns to a pair of vectors in $D_a$ the projection to
$T_a/D_a$ of the commutator of two local sections of $D$ that extend the
two vectors is called {\em Levi form}. A CR manifold is {\em Levi
nondegenerate} if the image of the Levi form coincides with
$T_a/D_a$ and if the null-space of the Levi form is trivial.

\item {\bf Uniformity.} A CR manifold is called {\em strongly
uniform} (cp. \cite{Miz}) if the Levi forms at different base
points are equivalent with respect to linear mappings.
\item {\bf Semi-simplicity of the structure group.} At any point
$a \in M$ the Levi form $\mathcal L_a$ can be used as a product of
a graded Lie algebra $\mathfrak g_-=\mathfrak g_{-1}\oplus
\mathfrak g_{-2}=D_a\oplus T_a/D_a$. The structure group $G$ is
the group of automorphisms of $\mathfrak g_-$ that preserve $J_a$
at $D_a$. For strongly uniform CR manifolds the structure groups
at different points are isomorphic.
\end{enumerate}

The nondegeneracy condition introduced above turns out to be too
strong and it excludes immediately CR manifolds whose codimension
$k$ is bigger than $n^2$.

In this paper we will study CR manifolds of codimension $k=2$ and
CR dimension $n=1$ which is the simplest case when $k>n^2$.
Instead of the nondegeneracy condition from above we will assume a
weaker condition introduced by V. Beloshapka in \cite{B2} that is based on 
the following observation: For given CR-dimension $n$ the dimension of the 
subspace of $T_aM$ that is generated by brackets of local sections of $D$ of order 
$\le r$ have precise upper bounds $k_{n,r}$. The generalized non-degeneracy 
requires that these upper bounds are attained until the complete tangent space is generated.

\begin{df}
Let $M$ be a CR-manifold of CR-dimension $n$ and codimension $k\ge n^2$. By $D^{(r)}_a$ 
denote the subspace of $T_aM$ generated by brackets of local sections of $D$ order $\le r$. 
Then $M$ is called non-degenerate at $a$ if there exists $r_0$ such that 
$T_a =D^{(r_0)}_a$ and $\dim D^{(r)}=k_{n,r}$ for $r < r_0$.
\end{df}

For a 4-dimensional manifold with a rank 2 distribution $D$ this
means that $TM$ is spanned by $D$ and commutators of first and
second orders of local sections of $D$. Such manifolds are called {\em
Engel manifolds} (cp. \cite{Eng}). In Section \ref{levt} we show
that Engel CR manifolds are strongly uniform in an appropriate
sense but their structure group is not semisimple. It follows that
the corresponding geometry is not parabolic.

Nevertheless, as the main result of this paper we obtain the
following

\begin{theo}
Engel CR manifolds admit a canonical Cartan connection.
\end{theo}

We prove this theorem in Section \ref{constr} by an explicit and
transparent construction. We will derive a complete system of
invariants that consists of four 1-dimensional components. A
geometric interpretation of these invariants will be given in
Section \ref{dist}.

In Section \ref{norm} we investigate the relation between the
invariants from the Cartan connection and a normal form earlier
constructed by the authors \cite{BES}.

\section{Engel manifolds and their Levi-Tanaka algebras}\label{levt} Let $M$
be an Engel manifold . Then we have a canonical flag of
distributions $D\subset D'\subset D''=TM$.

For two local sections $X,Y$ of $D$ defined in a neighbourhood of $a\in
M$ consider their bracket at $a$ followed the natural projection
$$\pi: T_aM \to T_aM/D_a.$$
The result takes values in $D'_a/D_a$ and depends only on $X(a)$
and $Y(a)$ (see Subsection \ref{brac}). Thus, we obtain an
algebraic bracket
$$\mathcal L_a^1:\, D_a\times D_a \to D'_a/D_a.$$

Analogously, we define a bilinear form
$$\tilde{\mathcal L}_a^1:\, D_a\times D_a' \to T_aM/D_a'$$
as the commutator of vector fields followed by the natural
projection. Since $\tilde{\mathcal L}_a^1(D_a,D_a)=0$ this bracket
lifts to a bracket
$$\mathcal L_a^2:\, D_a\times D_a'/D_a \to T_aM/D_a'.$$

For fixed $a$ denote $D_a$ by $\mathfrak g_{-1}$, $D_a'/D_a$ by
$\mathfrak g_{-2}$, and $T_aM/D_a'$ by $\mathfrak g_{-3}$. Then
$$\mathfrak g_-=\mathfrak g_{-1}\oplus\mathfrak
g_{-2}\oplus\mathfrak g_{-3}$$ forms a graded nilpotent Lie
algebra whose product is defined by $\mathcal L^1_a$ and $\mathcal
L^2_a$. This Lie algebra $\mathfrak g_-$ is called the {\em
Levi-Tanaka algebra} of $M$ at $a$.

Since $\mathfrak g_{-2}$ and $\mathfrak g_{-3}$ are 1-dimensional,
there is a unique direction $Y_a$ in $D_a=\mathfrak g_{-1}$ such
that $\mathcal L_a^2(Y_a,\mathfrak g_{-2})\equiv 0$. Thus, there
is a canonical line bundle $D^0$ in $D$ that consists of the
vectors that annihilate $\mathcal L_a^2$ (in the first argument).

It is well-known that Engel manifolds admit normal coordinates
$x,y,p,q$ such that
$$D= \Ann(dp-q\,dx,\; dy-p\,dx).$$

Then $D^0$ is spanned by $Y=\frac{\partial}{\partial q}$ and $D$
is spanned by $Y$ and one additional vector field of the form
$$X=\frac{\partial}{\partial x} +p\frac{\partial}{\partial y} +
q\frac{\partial}{\partial p} + B\frac{\partial}{\partial q},$$
where $B$ is a function of $x,y,p,q$. There is no canonical choice
of $X$ (i.e., of $B$)  from the Engel structure.

We remark that a choice of the direction field represented by the
complementary vector field $X$ on an Engel manifold is the same as
a third order ODE. In fact, the flow of $\xi$ consists of curves
$x(t),y(t),p(t),q(t)$ with $p=\frac{dy}{dx}$, $q=\frac{dp}{dx}$
and $\frac{dq}{dx}=B(x,y,p,q)$. Therefore,
$y'''(x)=B(x,y,y',y'')$.

The equivalence problem of third order ODE with respect to point
and contact transformations has been studied by E. Cartan
\cite{Carterc}, S.S. Chern \cite{Chern40} and H. Sato and A.Y.
Yoshikawa \cite{SY}.

Moreover, Engel CR manifolds are examples of multicontact sturctures that were introduced by M. Cowling, F. de Mari, A. Kor\'anyi and H.R. Reimann in \cite{CMKR}, see also \cite{Kor}.

Now let $M$ be an Engel CR manifold and let $Y$ be a smooth local
section of $D^0$. Then $X=JY$ represents a complementary direction field
to $Y$ in $D$. Note that the integrability condition is
automatically satisfied because the CR dimension is 1. Thus, the
CR structure encodes a third order ODE. Notice that the CR
structure cannot be completely recovered from the ODE since $J$
induces a scale of $X$ for any choice of a scale of $Y$.

The CR structure $J_a$ of an Engel CR manifold induces an
endomorphism on the $\mathfrak g_{-1}$-component of the
Levi-Tanaka algebra. This leads to a canonical family of bases
$(V_x,V_y,V_2,V_3)$ of $\mathfrak g_-$, where $V_y=tY$ annihilates
$\mathcal L^2$, and
\begin{align}\label{conJ}
V_x=&JV_y=tX\\\nonumber V_2=&[V_x,V_y]=-t^2[X,Y]\\\nonumber
V_3=&[V_x,V_2]=t^3[X,[X,Y]]\\\nonumber &[V_y,V_2]=0.
\end{align}
These conditions determine the basis up to $t\in \R^*$. It follows
that the Levi-Tanaka algebras of a Engel CR manifold at different
points are isomorphic.

\section{The cubic model}
The standard model of an Engel CR manifold is the 4-dimensional
real cubic $C$ in $\C^3$ defined in coordinates $z=x+\I
y,w_1=u_1+\I v_1,w_2=u_2+\I v_2$ by the equations
$$v_1=|z|^2, \qquad v_2=\Re z |z|^2.$$

We identify the cubic $C$ with a nilpotent Lie group $G_-$ of
holomorphic transformations on $\C^3$ which acts transitively on
$C$. This Lie group has a representation by matrices of the form
$$\begin{pmatrix} 1&z+\bar{z}&\I\bar{z}&2\I|z|^2+\I\bar{z}^2&w_2\\0&1&0&2\I\bar{z}
&w_1\\0&0&1&2z&z^2\\0&0&0&1&z\\0&0&0&0&1\end{pmatrix}
$$
where $\Im w_1=|z|^2$, $\Im w_2=\Re z|z|^2$.

The multiplication of two elements $(p,q_1,q_2)\in C $ and
$(z,w_1,w_2)\in C$ extends in the second factor to a holomorphic
action on $\C^3$ with coordinates $(z,w_1,w_2)$:
$$\begin{pmatrix} z\\ w_1 \\ w_2 \end{pmatrix}\mapsto
\begin{pmatrix} z+p\\ w_1 + q_1+ 2\I z\bar{p} \\ w_2 +q_2+ \I(2|p|^2+\bar{p}^2)z
+\I\bar{p} z^2+ (p+\bar{p})w_1 \end{pmatrix}. $$

The Lie algebra $\mathfrak g_-$ of $G_-$ coincides with the
Levi-Tanaka algebra of the homogeneous Engel CR manifold $C$. It
is represented by the matrix algebra of
$$\begin{pmatrix} 0&dz+d\bar{z}&\I d\bar{z}&0&du_2\\0&0&0&2\I d\bar{z}
&du_1\\0&0&0&2dz&0\\0&0&0&0&dz\\0&0&0&0&0\end{pmatrix}
$$
where $dz,d\bar{z}, du_1, du_2$ denote coordinates in the tangent
space $T_0 G_-$. We fix a basis that satisfies (\ref{conJ}):
\begin{align*}
V_x&=\begin{pmatrix} 0&1&\frac{\I}{2} &0&0\\0&0&0&\I
&0\\0&0&0&1&0\\0&0&0&0&\frac{1}{2}\\0&0&0&0&0\end{pmatrix} &\quad
V_y&=
\begin{pmatrix} 0&0 &-\frac{1}{2}&0&0\\0&0&0&-1 &0\\0&0&0&-\I &0\\0&0&0&0&-\frac{\I}{2}\\
0&0&0&0&0\end{pmatrix}\\
V_2&=\begin{pmatrix} 0&0&0 &0&0\\0&0&0&0
&1\\0&0&0&0&0\\0&0&0&0&0\\0&0&0&0&0\end{pmatrix} &\quad V_3&=
\begin{pmatrix} 0&0 &0&0&\frac{1}{2}\\0&0&0&0 &0\\0&0&0&0&0\\0&0&0&0&0\\
0&0&0&0&0\end{pmatrix}
\end{align*}
Or, if represented by vector fields,
\begin{align*}
V_x&=\frac{1}{2} \frac{\partial}{\partial x} +y
\frac{\partial}{\partial
u_1}+xy\frac{\partial}{\partial u_2}\\
V_y&=-\frac{1}{2} \frac{\partial}{\partial y}+x
\frac{\partial}{\partial
u_1}+\frac{3x^2+y^2}{2}\frac{\partial}{\partial u_2}\\
V_2&= \frac{\partial}{\partial u_1} +2x \frac{\partial}{\partial u_2} \\
V_3&=\frac{\partial}{\partial u_2}.
\end{align*}

The matrix Lie group $G_-$ can be extended to a bigger soluble Lie
group $G$ of holomorphic mappings that preserve $C$ by
incorporating the scaling
\begin{equation}\label{sca}
    (z,w_1,w_2)\mapsto (tz,t^2 w_1, t^3 w_2),
\end{equation}
which is the only remaining local automorphism of $C$ (see \cite{B1}). Thus, $G$
consists of matrices
$$g=\begin{pmatrix} t^3 &t^2 (p+\bar{p})&\I t^2 \bar{p}&\I t(2|p|^2+\bar{p}^2)&q_2\\
0&t^2&0&2\I t\bar{p}&q_1\\
0&0&t^2&2tp&p^2\\0&0&0&t&p\\0&0&0&0&1\end{pmatrix}
$$
The group $G$ is the automorphism group of the graded Lie algebra
$\mathfrak g_-$ with complex structure $J$ on $\mathfrak g_{-1}$.
Denote the subgroup of scalings in $G$ by $H$. Then the Lie
algebra $\mathfrak g$ of $G$ splits (as a vector space) into
$$\mathfrak g=\mathfrak g_- \oplus \mathfrak h$$
where $\mathfrak h$ is the one-dimensional Lie algebra of $H$. The
generator $V_0=t\frac{d}{dt}$ (in vector field representation) of
$\mathfrak h$ acts on $\mathfrak g_-$ as a grading element. Thus,
$\mathfrak g$ becomes a graded Lie algebra with $\mathfrak h=
\mathfrak g_0$.

The Maurer-Cartan form $\omega =g^{-1} dg$ of the Lie group $G$ in
matrix notation is
$$\begin{pmatrix} \frac{3dt}{t}& \frac{dz+ d\bar{z}}t&\frac{\I d\bar{z}}{t}&0&
\frac{du_2 - z du_1- \bar{z} du_1-\frac{\I}{2} z^2 d\bar{z}
+\frac{\I}{2} \bar{z}^2 dz}{t^3}\\
0& \frac{2dt}{t}&0&\frac{2\I d\bar{z}}{t} &\frac{du_1+\I z\,
d\bar{z}-\I \bar{z}\, dz}{t^2}
\\0&0& \frac{2dt}{t}& \frac{2\,dz}{t}&0
\\0&0&0&\frac{dt}{t}& \frac{dz}{t}\\0&0&0&0&0
\end{pmatrix}.$$

The forms
\begin{align*}
    \Phi_0&=\frac{dt}{t} \qquad \Phi_x=\frac{2dx}{t}  \qquad
    \Phi_y=-\frac{2dy}{t}\\
    \Phi_2&=\frac{du_1+\I z\, d\bar{z}-\I \bar{z}\, dz}{t^2}=
    \frac{du_1-2y\,dx+2x\,dy}{t^2}\\
    \Phi_3&=\frac{du_2 - z\, du_1- \bar{z}\, du_1-\frac{\I}{2} z^2 d\bar{z}
+\frac{\I}{2} \bar{z}^2 dz}{t^3}\\
&=\frac{du_2 -2x\,du_1+2xy\,dx - (x^2-y^2)dy}{t^3}
\end{align*}

constitute a dual basis to $(V_0,V_x,V_y,V_2,V_3)$. The structure
equation $d\omega+ \frac{1}{2}[\omega,\omega]=0$, or direct
verification, yields the structure formulae
\begin{align*}
d\Phi_0&=0, \qquad d\Phi_x= \Phi_x \wedge \Phi_0, \qquad
d\Phi_y=\Phi_y \wedge \Phi_0,\\
d\Phi_2&= 2\, \Phi_2\wedge \Phi_0-\Phi_x\wedge\Phi_y, \qquad
d\Phi_3= 3\, \Phi_3 \wedge \Phi_0- \Phi_x\wedge \Phi_2.
\end{align*}

\section{Preliminaries on Cartan connections}\label{prelim}
Our aim is to construct a canonical Cartan connection for Engel CR
manifolds based on the cubic model. This means that for a general
Engel CR manifold $M$ we have to find ``curved'' analogues of the
principal $\R^*$-bundle $G\to G/H=C$ and the $\mathfrak g$-valued
Maurer-Cartan form. More precisely, we will construct a principal
$\R^*$-bundle $\mathcal G\to M$ and a $\mathfrak g$-valued 1-form
$\omega:\, T\mathcal G \to \mathfrak g$ that satisfy the
properties

\begin{enumerate}
\item for any $p\in \mathcal G$ the mapping $\omega_p: T_p \mathcal G \to \mathfrak
g$ is isomorphic
\item the fundamental vector field of the principal action is
mapped to $V_0$
\item  $R^*_h \omega = \Ad(h^{-1}) \omega$ for any $h\in H$
\end{enumerate}

The curvature tensor of the Cartan connection is defined as
$$K=d\omega + \frac{1}{2}[\omega,\omega].$$
Differentiation of (iii) implies
$$\omega([X,Y])=[\omega(X),\omega(Y)]$$
when $X$ is the fundamental vector field of the principal action.
This means that $K$ is a horizontal form, i.e. it vanishes if one
of the arguments is vertical.

The curvature function on $\mathcal G$ is defined by
$$\kappa(p)(X,Y)=K(\omega|_p^{-1}X,\omega|_p^{-1}Y).$$
It that takes values in the 2-cochains $\mathcal C^2(\mathfrak
g_-,\mathfrak g)$. The $\mathfrak g$-module $\mathcal
C^2(\mathfrak g_-,\mathfrak g)$ splits into the subspaces of
homogeneity $i=-1,\dots,5$. Let $\partial$ be the cochain
operator. Then the $\partial\kappa$ can be expressed by the
so-called Bianchi-Identity (see \cite{CS})
$$ \partial \kappa(X,Y,Z)=-\sum_{cycl} \kappa(\kappa_-(X,Y),Z)+\omega^{-1}(X)\kappa(Y,Z),$$
where $X,Y,Z\in \mathfrak g_-$, $\kappa_-$ is the $\mathfrak
g_-$-part of $\kappa$ and the sum runs over all cyclic
permutations of $(X,Y,Z)$. It follows that the component $\partial
\kappa^{(i)}$ of homogeneity $i$ can be expressed by lower
components and therefore the components of degree $i$ are
determined by the lower components up to $\ker \partial=\mathcal
Z^2(\mathfrak g_-,\mathfrak g)$. Thus, the ``essential'' part of
the curvature takes values in $\mathcal Z^2(\mathfrak
g_-,\mathfrak g)$.

We will see below that there is a unique Cartan connection whose
$\partial$-exact part (with respect to some complementary
subspace) of the curvature function vanishes. Our choice of this
complementary subspace is based on the fact that the Lie algebra
$\mathfrak g$ splits into 1-dimensional subspaces. This induces a
splitting of $\mathcal C^2(\mathfrak g_-,\mathfrak g)$ into
1-dimensional components. We choose a basis for the complement
that involves a minimal number of those 1-dimensional components.
This leads to the simplest expressions for the curvatures (see
Subsection \ref{cohom} for the choice adopted in our
construction).

Since we have fixed a basis $(V_0,V_x,V_y,V_2,V_3)$ of $\mathfrak
g$, the connection form $\hat{\Phi}$ amounts to a distinguished
coframe
$\hat{\Phi}^0,\hat{\Phi}^x,\hat{\Phi}^y,\hat{\Phi}^2,\hat{\Phi}^3
$ on $\mathcal G$ or to its dual frame
$\hat{V}_0,\hat{V}_x,\hat{V}_y,\hat{V}_2,\hat{V}_3$. Notice that
$\hat{V}_0=\hat{\Phi}^{-1}(V_0)$,
$\hat{V}_x=\hat{\Phi}^{-1}(V_x)$,
$\hat{V}_y=\hat{\Phi}^{-1}(V_y)$,
$\hat{V}_2=\hat{\Phi}^{-1}(V_2)$,
$\hat{V}_3=\hat{\Phi}^{-1}(V_3)$.

The curvature components
$$K^\alpha_{\beta\gamma}:=K^\alpha(\hat{V}_\beta,\hat{V}_\gamma)$$
equal $-2k^\alpha_{\beta\gamma}$, with
$$k^\alpha_{\beta\gamma}=\hat{\Phi}^\alpha([\hat{V}_\beta,\hat{V}_\gamma]) -
\Phi^\alpha[V_\beta,V_\gamma]$$

where $\Phi^\alpha[V_\beta,V_\gamma]$ are the structure constants
of $\mathfrak g$. This shows that the curvature can be interpreted
as the aberration of the commutator relations between the frame
vector fields from the commutator relations between their images
under $\hat{\Phi}$ in $\mathfrak g$.

\section{Cartan Connection of Engel CR manifolds}\label{constr}

First, we define the Cartan bundle $\mathcal G:=D^0$, i.e.
$\mathcal G$ is the line bundle of all distinguished vectors $tY$.
It is clearly a $\R^*$-principal bundle. Notice that
$$\mathcal G\oplus J\mathcal G=D\subset D'\subset TM.$$

Now we construct a canonical frame field on $T\mathcal G$. First
of all we have the fundamental vector field $\hat{V}_0$ that is
induced from the principal action of $R^*$. Consider the two
mappings $\tilde{V}_x, \tilde{V}_y :\, \mathcal G\to TM$, defined
by
$$\tilde{V}_x(Y)=JY,\quad \tilde{V}_y(Y)=Y$$
If a connection form $\hat{\Phi}^0$ at $\mathcal G$ with
$\hat{\Phi}^0(\hat{V}_0)=1$ is given, these mappings lift to
vector fields on $\mathcal G$
\begin{align*}
\hat{V}_x&=\tilde{V}_x + a_{x0} \hat{V}_0\\
\hat{V}_y&=\tilde{V}_y + a_{y0} \hat{V}_0
\end{align*}
by the condition $\hat{\Phi}^0(\hat{V}_\alpha)=0$. We try to
determine $\hat{\Phi}^0$ and the missing vector fields
$\hat{V}_2,\hat{V}_3$ by making the commutator bracket relations
as ``close as possible'' to the corresponding relations for the
cubic model, i.e. we make as many curvature components vanishing
as possible.

For calculations we fix a local section $Y$ on $\mathcal G$ and
the frame $T_x=JY$, $T_y=Y$, $T_2=[T_x,T_y]$, $T_3=[T_x, T_2]$ on
$M$. Denote by $\varphi^x,\varphi^y,\varphi^2,\varphi^3$ the dual
coframe. We may choose the scale of $T_y$ in such a way that
$$[T_y,T_2]\in D.$$
In fact,
$$[\tau T_y,[\tau T_x,\tau T_y]]\equiv\tau^3[T_y,T_2]+
3\tau^2(T_y\tau)T_2\quad \mod D.$$ Thus, $\tau$ has to satisfy a
differential equation $T_y(\log
\tau)=\frac{1}{3}\varphi^2([T_y,T_2])$.

This choice and the Jacobi identity imply
\begin{align*}
[T_x,T_y]=&T_2\\
[T_x,T_2]=&T_3\\
[T_y,T_2]=&\phi_{y2}^xT_x+\phi_{y2}^yT_y\\
[T_x,T_3]=&\phi_{x3}^xT_x+\phi_{x3}^yT_y+\phi_{x3}^2T_2
+\phi_{x3}^3T_3\\
[T_y,T_3]=&[T_x,[T_y,T_2]]=(T_x\phi_{y2}^x)T_x+
(T_x\phi_{y2}^y)T_y+\phi_{y2}^yT_2\\
[T_2,T_3]=&[[T_x,T_3],T_y]+[T_x,[T_y,T_3]]=\\
&(-T_y\phi_{x3}^x+T^2_x\phi_{y2}^x-\phi_{y2}^x\phi_{x3}^2-
\phi_{x3}^3T_x\phi_{y2}^x)T_x+\\
&+(-T_y\phi_{x3}^y+T^2_x\phi_{y2}^y-\phi_{y2}^y\phi_{x3}^2-
\phi_{x3}^3T_x\phi_{y2}^y)T_y+\\
&+(\phi_{x3}^x+2T_x\phi_{y2}^y-T_y\phi_{x3}^2-
\phi_{x3}^3\phi_{y2}^y)T_2+\\
&+(\phi_{y2}^y-T_y\phi_{x3}^3)T_3
\end{align*}

Here
$\phi^\alpha_{\beta\gamma}=\varphi^\alpha([T_\beta,T_\gamma])$.
Further consequences of the Jacobi identity which will be used in
calculations are
\begin{align}\label{jacobi}
0=&[T_x,[T_2,T_3]]+[T_2,[T_3,T_x]]\\\nonumber
0=&[T_y,[T_2,T_3]]+[T_2,[T_3,T_y]]+[T_3,[T_y,T_2]].
\end{align}

By $t$ we denote the fibre coordinate in $\mathcal G$ with respect
to the trivialization induced by the fixed section $Y$, i.e. for
$W\in \mathcal G$ the coordinate $t(W)=\varphi^y(W)$. Then
\begin{align*}
\hat{V}_0&=t\frac{\partial}{\partial t}\\
\hat{V}_x&=tT_x + a_{x0} \hat{V}_0\\
\hat{V}_y&=tT_y + a_{y0} \hat{V}_0\\
\hat{V}_2&=a_{22}T_2 + a_{2x} T_x +a_{2y}
T_y + a_{20}\hat{V}_0\\
\hat{V}_3&=a_{33}T_3 + a_{32}T_2 + a_{3x} T_x +a_{3y} T_y +
a_{30}\hat{V}_0
\end{align*}
The dual coframe is then
\begin{align*}
\hat{\Phi}^0&=\frac{dt}{t}+ b_{0x}\varphi_x + b_{0y}\varphi_y + b_{02}\varphi_2 + b_{03}\varphi_3\\
\hat{\Phi}^x&= \frac{1}{t}\varphi_x +b_{x2}\varphi_2+b_{x3}\varphi_3 \\
\hat{\Phi}^y&= \frac{1}{t}\varphi_y +b_{y2}\varphi_2+b_{y3}\varphi_3 \\
\hat{\Phi}^2&=b_{22}\varphi_2 + b_{23} \varphi_3 \\
\hat{\Phi}^3&=b_{33}\varphi_{3}
\end{align*}
with $b_{33}=\frac{1}{a_{33}}$, $b_{22}=\frac{1}{a_{22}}$, etc.

The invariance condition (iii) of a Cartan connection means
$$[\hat{V}_0,\hat{V}_j]=-|j|\,\hat{V}_j.$$ Here $|j|=j$ for $j=0,2,3$ and
$|j|=1$ for $j=x,y$. It follows that $a_{j*}(m,t)= t^{|j|}
\alpha_{j*}(m)$ and $b_{j*}(m,t)= t^{-|j|} \beta_{j*}(m)$. In
particular,
$$\hat{\Phi}_0=\frac{dt}{t} + \varpi$$
where $\varpi$ is a 1-form on $M$. (More precisely, $\varpi$ is
the pullback of a 1-form at $M$.)

Now we will investigate all curvature components. By our
construction all components of homogeneity $-1$ and $0$ vanish
automatically:
\begin{align*}
k_{xy}^3=\hat{\Phi}^3([\hat{V}_x,\hat{V}_y])=0\\
k_{y2}^3=\hat{\Phi}^3([\hat{V}_y,\hat{V}_2])=0
\end{align*}
By setting $a_{22}=t^2$ and $a_{33}=t^3$ we annihilate
\begin{align*}
k_{xy}^2=\hat{\Phi}^2([\hat{V}_x,\hat{V}_y]-\hat{V}_2)&=
b_{22}(t^2-a_{22})\\
k_{x2}^3=\hat{\Phi}^3([\hat{V}_x,\hat{V}_2]-\hat{V}_3)&=
b_{33}(t^3-a_{33}).
\end{align*}
Next we consider the components of homogeneity $1$. The subspace
of closed cochains is 5-dimensional and coincides with the
subspace of coboundaries (see Subsection \ref{cohom} for the
explicit calculation of the cohomologies). Thus vanishing of
curvature of homogeneity $1$ yields
\begin{align*}
k_{xy}^y=\hat{\Phi}^y([\hat{V}_x,\hat{V}_y]-\hat{V}_2)&=
a_{x0}-\frac{a_{2y}}{t}\\
k_{xy}^x=\hat{\Phi}^x([\hat{V}_x,\hat{V}_y]-\hat{V}_2)&=
-a_{y0}-\frac{a_{2x}}{t}\\
\hat{\Phi}^2([\hat{V}_x,\hat{V}_2]-\hat{V}_3)&=2 a_{x0}+
\frac{a_{2y}}{t}-\frac{a_{32}}{t^2}\\
\hat{\Phi}^2([\hat{V}_y,\hat{V}_2])&=2 a_{y0}- \frac{a_{2x}}{t}\\
\hat{\Phi}^3([\hat{V}_x,\hat{V}_3])&=3a_{x0}+
\frac{a_{32}}{t^2}+t\phi^3_{x3}\\
\hat{\Phi}^3([\hat{V}_y,\hat{V}_3])&=3a_{y0}
\end{align*}
We find
\begin{align*}
a_{32}&=-\frac{t^3}{2}\phi^3_{x3},\\
a_{2y}&=-\frac{t^2}{6}\phi^3_{x3},\\
a_{x0}&=-\frac{t}{6}\phi^3_{x3},\\
a_{2x}&=0, \\
a_{y0}&=0
\end{align*}

Due to the Bianchi identity the curvature of homogeneity 2 takes
values in $\mathcal Z^{(2)}$ which has dimension 6. It splits into
an exact and a non-exact component. Below, we used a bar to
separate the exact part at the left hand sind from the non-exact
part at the right hand side.
\begin{align*}
k^0_{xy} &= \left.\frac{1}{6}t^2T_y \phi^3_{x3}-a_{20}\right|\\
k^x_{x2} &= \left.-a_{20}-\frac{1}{t}a_{3x}\right|\\
k^y_{x2} &= \left.\frac{1}{18}t^2(\phi^3_{x3})^2-
\frac{1}{6}t^2T_x\phi^3_{x3}-\frac{1}{t}a_{3y}\right|\\
k^x_{y2} &= \bigg|t^2\phi^x_{y2}\\
k^y_{y2} &= \left.\frac{1}{6}t^2T_y\phi^3_{x3}-a_{20}\right|-
\frac{1}{3}t^2T_y\phi^3_{x3}+t^2\phi^y_{y2}\\
k^2_{x3} &=\frac{1}{t}a_{3y}+\frac{1}{6}t^2T_x\phi^3_{x3}
-\frac{1}{18}t^2(\phi^3_{x3})^2\left|+\frac{11}{36}t^2(\phi^3_{x3})^2-
\frac{2}{3}t^2T_x\phi^3_{x3}+t^2\phi^2_{x3}\right.\\
k^2_{y3}&=-\frac{1}{t}a_{3x}-\frac{1}{6}t^2T_y\phi^3_{x3}\left|-
\frac{1}{3}t^2T_y\phi^3_{x3}+t^2\phi^y_{y2}\right.\\
k^3_{23}&= 3a_{20}-\frac{1}{t}a_{3x}-\frac{2}{3}t^2T_y\phi^3_{x3}
\left|-
\frac{1}{3}t^2T_y\phi^3_{x3}+t^2\phi^y_{y2}\right.\\
\end{align*}

The functions $a_{3x}$ and $a_{3y}$ and $a_{20}$ will be
determined by equating the exact part of homogeneity 2 to zero
\begin{align*}
a_{3x}& =-\frac{t^3}{6}T_y\phi^3_{x3}\\
a_{3y}&=\frac{t^3}{18}(\phi^3_{x3})^2-\frac{t^3}{6}T_x\phi^3_{x3}\\
a_{20} &= \frac{t^2}{6}T_y\phi^3_{x3}.
\end{align*}

The essential curvatures of homogeneity 2 are
\begin{align*}
R^x_{y2}=&t^2\phi^x_{y2}\\
R^y_{y2}=&t^2\phi^y_{y2} -\frac{t^2}{3}T_y\phi^3_{x3}\\
R^2_{x3}=&t^2\phi^2_{x3}+\frac{11t^2}{36}(\phi^3_{x3})^2-\frac{2t^2}{3}T_x\phi^3_{x3}
\end{align*}

The curvature component of homogeneity 3 splits into a
1-dimensional exact part, a 1-dimensional non-exact closed part
and a 5-dimensional non-closed part. Below, we used two bars to
separate the three components.
\begin{align*}
k^0_{x2} =& B_3\bigg|\bigg|\\
k^0_{y2}=& \left|\left| \frac{1}{4}\hat{V}_x
R^x_{y2}+\frac{1}{4}\hat{V}_y R^y_{y2}
\right.\right.\\
k^x_{x3}=&  B_3\left|\left|
-\frac{3}{2}\hat{V}_x R^y_{y2}+\frac{1}{2}\hat{V}_y R^2_{x3}\right.\right.\\
k^y_{x3} =&\bigg| R^y_{x3}\bigg|\\
k^x_{y3} =& \left|\left| \frac{3}{4}\hat{V}_x
R^x_{y2}-\frac{1}{4}\hat{V}_y R^y_{y2}
\right.\right.\\
k^y_{y3} =&B_3\bigg|\bigg|+\hat{V}_x R^y_{y2}\\
k^2_{23}=&2B_3\left|\left|+\frac{1}{2}\hat{V}_x
R^y_{y2}-\frac{1}{2}\hat{V}_y R^2_{x3}\right.\right. ,
\end{align*}
with
$$B_3:=-a_{30}+\frac{1}{3}t^3T_xT_y\phi^3_{x3}
-\frac{1}{6}t^3T_y T_x\phi^3_{x3}-
\frac{1}{18}t^3\phi^3_{x3}T_y\phi^3_{x3}.$$

The only essential curvature of homogeneity 3 is
$$R^y_{x3}=\phi^y_{x3}t^3+\frac{t^3}{12}\phi^3_{x3}T_x\phi^3_{x3}-
\frac{t^3}{6}T_x^2\phi^3_{x3}+ \frac{5t^3}{216} (\phi^3_{x3})^3+
\frac{t^3}{6}\phi^3_{x3}\phi^2_{x3}$$

Vanishing of the exact part yields
$$a_{30}=\frac{t^3}{3}T_xT_y\phi^3_{x3}
-\frac{t^3}{6}T_y T_x\phi^3_{x3}-
\frac{t^3}{18}\phi^3_{x3}T_y\phi^3_{x3}.$$

In homogeneity 4 and 5 the cochain operator $\partial$ has no
kernel. Therefore, according to the Bianchi identity, all
curvature of these homogeneities can be expressed in terms of
curvature of lower homogeneity:
\begin{align*}
k^0_{x3} =&
-\frac{1}{2}\hat{V}_x^2R^y_{y2}-\frac{1}{3}\hat{V}_y\hat{V}_x
R^2_{x3}+\frac{1}{2} \hat{V}_x\hat{V}_y
R^2_{x3}+\frac{1}{3}\hat{V}_y R^y_{x3}\\
k^0_{y3} =& \frac{1}{4}\hat{V}_x^2
R^x_{y2}+\frac{1}{4}\hat{V}_x\hat{V}_y R^y_{y2}\\
k^x_{23} =& -\frac{1}{2} \hat{V}_x\hat{V}_y R^y_{y2}+\frac{1}{2}
\hat{V}_x^2 R^x_{y2}+ \frac{3}{2} \hat{V}_y\hat{V}_x
R^y_{y2}-\frac{1}{2} \hat{V}_y^2 R^2_{x3}
-R^x_{y2}R^2_{x3}\\
k^y_{23} =& \frac{1}{2} \hat{V}_x^2 R^y_{y2}- \frac{1}{3}
\hat{V}_y\hat{V}_x R^2_{x3}+\frac{1}{2} \hat{V}_x\hat{V}_y
R^2_{x3}-\frac{2}{3} \hat{V}_y R^y_{x3}-R^y_{y2}R^2_{x3}\\
k^0_{23}=&\frac{1}{2}\hat{V}_x^3 R^x_{y2}+3\hat{V}_x
\hat{V}_y\hat{V}_x R^y_{y2}-\frac{3}{2} \hat{V}_y\hat{V}_x^2
R^y_{y2}+ \frac{1}{2}\hat{V}_y\hat{V}_x\hat{V}_y R^2_{x3}\\&
-\hat{V}_x\hat{V}_y^2 R^2_{x3}-\frac{1}{2}\hat{V}_x^2\hat{V}_y
R^y_{y2}
 -R^x_{y2} \hat{V}_x R^2_{x3}-  R^2_{x3}\hat{V}_x R^x_{y2}\\&+\frac{1}{2} R^y_{y2} \hat{V}_y R^2_{x3}
 -\frac{3}{2} R^y_{y2} \hat{V}_x R^y_{y2}+ R^x_{y2} R^y_{x3}
\end{align*}

Here we used the Jacobi identity (\ref{jacobi}). We conclude

\begin{cor}\label{cor1}
An Engel CR manifold is locally equivalent to the cubic $C$ if and
only if the invariants $R^x_{y2}, R^y_{y2},R^2_{x3},R^y_{x3}$
vanish identically.
\end{cor}

The structure equations take the form
\begin{align*}
d\hat{\Phi}^3&=3 \hat{\Phi}^3\wedge \hat{\Phi}^0
-\hat{\Phi}^x\wedge \hat{\Phi}^2+K^3\\
d\hat{\Phi}^2&=2 \hat{\Phi}^2\wedge \hat{\Phi}^0
-\hat{\Phi}^x\wedge \hat{\Phi}^y+K^2\\
d\hat{\Phi}^y&=\hat{\Phi}^y\wedge \hat{\Phi}^0+K^y\\
d\hat{\Phi}^x&=\hat{\Phi}^x\wedge \hat{\Phi}^0+K^x\\
d\hat{\Phi}^0&=K^0=\pi^*(d\varpi)
\end{align*}

The curvature form $K^0$ is a pullback of a 2-form on $M$. Its
vanishing is equivalent to the integrability of the horizontal
distribution defined as the null distribution of the connection
form
$$\hat{\Phi}^0=\frac{dt}{t} + \frac{\phi^3_{x3}}{6}\varphi_x-\frac{T_y\phi^3_{x3}}{6}\varphi_2-
\frac{2T_xT_y\phi^3_{x3}-T_yT_x\phi^3_{x3}}{6}\varphi_3.$$

\section{Distinguished frames and interpretation of the curvature} \label{dist}
The Cartan connection constructed above induces a family of
distinguished frames at any point $a$ of the underlying manifold
itself. Let $\pi: \mathcal G \to M$ be the projection of the
Cartan bundle and $\hat{a}\in \pi^{-1}(a)$. Then for any choice of
$\hat{a}$ we define for $\alpha=x,y,2,3$
$$\tilde{V}_\alpha=d\pi|_{\hat{a}} \hat{V}_\alpha$$
Another choice of $\hat{a}'=t\hat{a}$ leads to the frame
$$\tilde{V}_x'=t\tilde{V}_x, \quad \tilde{V}_y'=t\tilde{V}_y,\quad \tilde{V}_2'=t^2\tilde{V}_2,
\quad \tilde{V}_3'=t^3\tilde{V}_3.$$ Thus, the Cartan connection
induces a complete splitting of the tangent spaces of $M$. The
essential curvatures $R^x_{y2}, R^y_{y2}, R^y_{x3}, R^2_{x3},$ can
be interpreted as the $\tilde{V}_x$ component of
$[\tilde{V}_y,\tilde{V}_2]$, the $\tilde{V}_2$ component of
$[\tilde{V}_y,\tilde{V}_3]$ and the $\tilde{V}_y$ and
$\tilde{V}_2$ components of $[\tilde{V}_x,\tilde{V}_3]$,
respectively. This proves the following geometric interpretation
of the curvatures

\begin{cor} Vanishing of $R^x_{y2}$ is equivalent to the integrability
of the distribution spanned by $\tilde{V}_y,\tilde{V}_2$,
vanishing of $R^y_{x3}$ and $R^2_{x3}$  is equivalent to the
integrability of the distribution spanned by
$\tilde{V}_x,\tilde{V}_3$. If $\hat{V}_xR^x_{y2}\equiv 0$ then
$R^y_{y2}\equiv 0$ is equivalent to the integrability of the
distribution spanned by $\tilde{V}_y,\tilde{V}_3$.
\end{cor}

The splitting of $TM$ can be understood in terms of $M$ only in
the following special case: Choose $T_x,T_y$ as above. The
remaining freedom is a rescaling of $T_x,T_y$ by a multiplier $f$
with $T_y f=0$. Suppose there is a solution of
\begin{equation}\label{dex}
\frac{Xf}{f}=X(\log f)=\frac{\phi^3_{x3}}{3}. \end{equation} This
determines all partial derivatives of $f$, thus $f$ is determined
up to a multiplicative constant. In this case $\tilde{V}_x=tfT_x$,
$\tilde{V}_y=tfT_y$, $\tilde{V}_2=[\tilde{V}_x,\tilde{V}_y]$,
$\tilde{V}_3=[\tilde{V}_x,\tilde{V}_2]$. All formulae simplify
significantly. In particular $\hat{\Phi}^0=\frac{dt}{t}$ and
$d\hat{\Phi}^0=0$. Vice versa, let $d\hat{\Phi}^0=0$. Then
$d\varpi=0$ and locally $\varpi=d\psi$ for some function $f$. But
then $T_y\psi=0$ and $T_x\psi=\frac{\phi^3_{x3}}{6}$ and, hence,
$f=\exp(2\psi)$ solves equation (\ref{dex}). It follows

\begin{cor}
On an Engel CR manifold M there exists locally a unique family of
vector fields $tV_x,tV_y,t^2V_2,t^3V_3$ ($t\in\R^*$) on $M$ such
that
\begin{align*}
[V_x,V_y]=&V_2\\
[V_x,V_2]=&V_3\\
[V_y,V_2]\in& D\\
[V_x,V_3]\in& D'
\end{align*}
if and only if $d\varpi=0$ which is equivalent to $K^0=0$.
\end{cor}

Obviously, the families $(tV_x, tV_y, tV_2, tV_3)$ and
$(t\tilde{V}_x, t\tilde{V}_y, t\tilde{V}_2, t\tilde{V}_3)$
coincide.

\section{Relation to the normal form}\label{norm}
In \cite{BES} the authors constructed a normal form for embedded
real-analytic Engel CR manifolds. For any such manifold $M$ there
exist local coordinates $(z,w_1,w_2)$ in a neighbourhood of a
point $a\in M$ in the ambient space such that the equation of $M$
takes the form
\begin{align*}
\Im w_1=&|z|^2+A_1\Re z^2\bar{z}^3+A_2 \Im z^2\bar{z}^3+\cdots \\
\Im w_2=&\Re z^2\bar{z}+B_1\Re z^4\bar{z}+ B_2 \Re z^2\bar{z}^3+
B_3 \Im z^2\bar{z}^3 +B_4 \Re z^5\bar{z}\\& B_5 \Im z^5\bar{z}+
B_6 \Re z^4\bar{z}^2 +B_7 \Im z^4\bar{z}^2+B_8|z|^6 +\cdots
\end{align*}
where the dots indicate terms of higher homogeneity.

The expressions of the invariants in terms of the normal form are
\begin{align*}
R^x_{y2}(0)&=(2B_1-B_2)t^2\\
R^y_{y2}(0)&=-3B_3t^2\\
R^2_{x3}(0)&=(2B_1-5B_2)t^2\\
R^y_{x3}(0)&=(4A_1+5B_4-2B_6-6B_8)t^3
\end{align*}

The distinguished frames defined in Section \ref{dist} at the
reference point $a$ are
$$\frac{t}{2}\frac{\partial}{\partial x},\quad -\frac{t}{2}\frac{\partial}{\partial y},
\quad t^2\frac{\partial}{\partial u_2}, \quad
t^3\frac{\partial}{\partial u_3} -t^3(B_1-B_2)
\frac{\partial}{\partial y}.$$ Notice, that the normal coordinates
themselves provide a family of distinguished frames at the
reference point, namely (for the most natural choice)
$$t\frac{\partial}{\partial x},\quad t\frac{\partial}{\partial y},
\quad t^2\frac{\partial}{\partial u_2}, \quad
t^3\frac{\partial}{\partial u_3}.$$ Since $B_1-B_2$ is an
invariant of homogeneity 2 the frame from the Cartan connection
and the frame from the normal form are equivalent.

It is natural to introduce the notion of umbilicity. In analogy to the definition of umbilic points of hypersurface given in \cite{CM} we call a point $a$ of an Engel CR-manifold umbilic if $R^x_{y2}(a)= R^y_{y2}(a)=R^2_{x3}(a)=R^y_{x3}(a)=0$. Then Corollary \ref{cor1} can be reformulated as \medskip

{\bf Corollary \ref{cor1}'.} {\em 
An Engel CR manifold is locally equivalent to the cubic $C$ if and
only if all its points are umbilic.} \medskip

In terms of the normal form umbilicity is equivalent to vanishing of $B_1,B_2,B_3$ and $4A_1+5B_4-2B_6-6B_8$.

\section{Appendix}
\subsection{Algebraic brackets}\label{brac}
Let $D\subset TM$ be a distribution of rank $n$ on a
$n+k$-dimensional manifold $M$ and let $\mathcal L_a$ be the
bilinear mapping
$$\mathcal L_a: (X,Y) \to \pi([X,Y]_a)$$
where $X,Y$ are local sections of $D$, $a\in M$ and $\pi$ is the natural
projection $T_aM \to T_aM/D_a$. Then $\mathcal L_x(X,Y)$ depends
only on $X_a$ and $Y_a$. Indeed, let $\theta$ be an $\R^k$ valued
1-form with $\ker \theta=D$. Then $\theta$ identifies $T_aM/D_a$
with $\R^k$ and, thus $\mathcal L_a$ with
$$\theta_a([X,Y])=-2d\theta_a(X,Y)+X\theta(Y)-Y\theta(X).$$
Since $\theta(X)=\theta(Y)=0$, we conclude
$\theta_a([X,Y])=-2d\theta_a(X_a,Y_a)$.

\subsection{The space of 2-cocycles}\label{cohom}

Let  $V_a$ a basis of $\mathfrak g$ and $\Phi^a$ the dual basis,
where $a$ runs over $\{0,x,y,2,3\}$. Then a 2-cochain can be
written as $\Phi=\sum \phi_{\alpha\beta}^a \,
\Phi^\alpha\wedge\Phi^\beta\otimes V_a$ (the indices
$\alpha,\beta$ run over $\{x,y,2,3\}$ and
$\phi_{\alpha\beta}^a=-\phi_{\beta\alpha}^a$). Thus, the space of
2-cochains $\mathcal C^2(\mathfrak g_-,\mathfrak g)$ has dimension
30. The subspace $\mathcal Z^2$ of $\partial$-closed cochains is
17-dimensional and is characterized by the following conditions
(in order of their homogeneity):
\begin{align*}
\phi_{xy}^x -\phi_{y2}^2+\phi_{y3}^3&=0 &&&&\\
\phi_{x2}^x +\phi_{y2}^y +5\phi_{xy}^0-\phi_{23}^3&=0 &
\phi_{y3}^2+ 3\phi_{xy}^0-\phi_{23}^3&=0 &&\\
\phi_{y3}^y- \phi_{x3}^x&=0 & \phi_{y2}^0&= 0 &\phi_{y3}^x&= 0 \\
\phi_{23}^2- 2\phi_{x3}^x&=0 & \phi_{x2}^0+ \phi_{x3}^x&=0 &&\\
\phi_{y3}^0= 0 \qquad \phi_{23}^x&= 0 &
 \phi_{x3}^0 &= 0 & \phi_{23}^y&= 0\\\phi_{23}^0 &= 0 &&
\end{align*}

The subspace $\mathcal B^2$ of $\partial$-exact cochains is
described by the additional equations
\begin{align*}
\phi^x_{y2}&=0 \quad \phi^y_{x2}+\phi^2_{x3}=0 \quad \phi^y_{y2}- \phi^0_{xy}=0\\
\phi^y_{x3}&=0
\end{align*}

Thus, the cohomology $H^2(\mathfrak g_-,\mathfrak g)$ is
4-dimensional and can be represented by the following
complementary subspaces to $\mathcal B^2(\mathfrak g_-,\mathfrak
g)$ in $\mathcal Z^2(\mathfrak g_-,\mathfrak g)$: in homogeneity 2
we choose the subspace spanned by
\begin{align*}
&\Phi^y\wedge \Phi^2\otimes V_x, \qquad \Phi^x\wedge \Phi^3\otimes V_2\\
&\Phi^y\wedge \Phi^2\otimes V_y + \Phi^y\wedge \Phi^3\otimes V_2
+\Phi^2\wedge \Phi^3\otimes V_3,
\end{align*}
in homogeneity 3 we choose the subspace spanned by
$$\Phi^x\wedge \Phi^3\otimes V_y.$$

\subsection{The curvatures}

Here we list the curvatures by their homogeneity
\begin{align*}
(-1):&\quad&  k_{xy}^3=& 0 & & & &\\
(0):&\quad& k_{y2}^3=&0, & k_{xy}^2=&0, &  k_{x2}^3=&0\\
(1):&\quad& k_{xy}^y=&0, & k_{xy}^x=&0, & k_{x2}^2=&0,\\
&&k_{y2}^2=&0,& k_{x3}^3=&0, & k_{y3}^3=&0\\
(2):&\quad& k^0_{xy} =&0, & k^x_{x2} =& 0, & k^y_{x2} =&0 \\
&& k^x_{y2} = & R^x_{y2}, & k^y_{y2} = &R^y_{y2}, & k^2_{x3} = &R^2_{x3}\\
&&k^2_{y3}=&R^y_{y2}, & k^3_{23}=&  R^y_{y2} & & \\
(3):&\quad&k^0_{x2} =&0, & k^0_{y2}=& \frac{1}{4}\hat{V}_x
R^x_{y2}+\frac{1}{4}\hat{V}_y R^y_{y2}, & k^x_{x3}=&
-\frac{3}{2}\hat{V}_x R^y_{y2}+\frac{1}{2}\hat{V}_y
R^2_{x3}\\
 & &k^y_{x3} =& R^y_{x3}, & k^x_{y3} =& \frac{3}{4}\hat{V}_x
R^x_{y2}-\frac{1}{4}\hat{V}_y R^y_{y2}, & k^y_{y3} =&\hat{V}_x R^y_{y2}\\
&&&&k^2_{23}=&\frac{1}{2}\hat{V}_x R^y_{y2}-\frac{1}{2}\hat{V}_y
R^2_{x3} &&
\end{align*}
\begin{align*}
(4):& \quad& k^0_{x3} =&
-\frac{1}{2}\hat{V}_x^2R^y_{y2}-\frac{1}{3}\hat{V}_y\hat{V}_x
R^2_{x3}+\frac{1}{2} \hat{V}_x\hat{V}_y
R^2_{x3}+\frac{1}{3}\hat{V}_y R^y_{x3}\\
&& k^0_{y3} =& \frac{1}{4}\hat{V}_x^2
R^x_{y2}+\frac{1}{4}\hat{V}_x\hat{V}_y R^y_{y2}\\
&& k^x_{23} =& -\frac{1}{2} \hat{V}_x\hat{V}_y
R^y_{y2}+\frac{1}{2} \hat{V}_x^2 R^x_{y2}+ \frac{3}{2}
\hat{V}_y\hat{V}_x R^y_{y2}-\frac{1}{2} \hat{V}_y^2 R^2_{x3}
-R^x_{y2}R^2_{x3}\\
&& k^y_{23} =& \frac{1}{2} \hat{V}_x^2 R^y_{y2}- \frac{1}{3}
\hat{V}_y\hat{V}_x R^2_{x3}+\frac{1}{2} \hat{V}_x\hat{V}_y
R^2_{x3}-\frac{2}{3} \hat{V}_y R^y_{x3}-R^y_{y2}R^2_{x3}
\end{align*}
\begin{multline*}
(5):\quad k^0_{23}=\frac{1}{2}\hat{V}_x^3 R^x_{y2}+3\hat{V}_x
\hat{V}_y\hat{V}_x R^y_{y2}-\frac{3}{2} \hat{V}_y\hat{V}_x^2
R^y_{y2}+
\frac{1}{2}\hat{V}_y\hat{V}_x\hat{V}_y R^2_{x3}-\\
-\hat{V}_x\hat{V}_y^2 R^2_{x3}-\frac{1}{2}\hat{V}_x^2\hat{V}_y
R^y_{y2}
 -R^x_{y2} \hat{V}_x R^2_{x3}-\\-  R^2_{x3}\hat{V}_x R^x_{y2}+\frac{1}{2} R^y_{y2} \hat{V}_y R^2_{x3}
 -\frac{3}{2} R^y_{y2} \hat{V}_x R^y_{y2}+ R^x_{y2} R^y_{x3}
\end{multline*}

The essential curvatures evaluate as
\begin{align*}
R^x_{y2}=&t^2\phi^x_{y2}\\
R^y_{y2}=&t^2\phi^y_{y2} -\frac{t^2}{3}T_y\phi^3_{x3}\\
R^2_{x3}=&t^2\phi^2_{x3}+\frac{11t^2}{36}(\phi^3_{x3})^2-\frac{2t^2}{3}T_x\phi^3_{x3}\\
R^y_{x3}=&t^3\phi^y_{x3}+\frac{t^3}{12}\phi^3_{x3}T_x\phi^3_{x3}-
\frac{t^3}{6}T_x^2\phi^3_{x3}+ \frac{5t^3}{216} (\phi^3_{x3})^3+
\frac{t^3}{6}\phi^3_{x3}\phi^2_{x3}
\end{align*}

%\subsection{Auxiliary vector fields in normal form}
%
%The vector fields $T_x,T_y,T_2,T_3$ in normal form are (up to 2nd
%order):
%\begin{align*}
%T_x&=\frac{1}{2} \frac{\partial}{\partial x} +y
%\frac{\partial}{\partial
%u_1}+xy\frac{\partial}{\partial u_2}\\
%T_y&=-\frac{1}{2} \frac{\partial}{\partial y}+x
%\frac{\partial}{\partial
%u_1}+\frac{3x^2+y^2}{2}\frac{\partial}{\partial u_2}\\
%T_2&= \frac{\partial}{\partial u_1} +2x \frac{\partial}{\partial u_2} \\
%T_3&=\frac{\partial}{\partial u_2}.
%\end{align*}
\bigskip

{\bf Acknowledgement.} The authors thank Keizo Yamaguchi and
Andrea Spiro for helpful discussions. Parts of this article were
written while the third named author stayed at Stockholm
University. He gratefully acknowledges the hospitality and
interesting discussions with Mikael Passare and August Tsikh.

\end{document}